\documentclass[reqno]{amsart}
\usepackage{amsmath,amsfonts,amssymb,amsthm}
\usepackage{hyperref}
\usepackage{mystyle}  

\usepackage{dsfont}
\def\1{{\mathds 1}}
\def\Z{{\mathds Z}}
\def\k{{\mathds k}}

\numberwithin{Lem}{section}
\numberwithin{Prop}{section}

\begin{document}

\bibliographystyle{plain}

\author{Michael Roitman}

\address{University of Michigan\\
Department of Mathematics\\
Ann Arbor, MI 48109-1109}
\email{roitman@umich.edu}

\title{On embedding of Lie conformal algebras into associative
  conformal algebras}

\date\today




\maketitle
\section*{Introduction}
\subsection*{Conformal algebras}
A conformal algebra is, roughly speaking, a linear space $\goth A$ with 
infinitely many bilinear products $(n):\goth A\times \goth A\to \goth
A$, parameterized
by an non-negative integer $n$, and a derivation $D:\goth A\to \goth A$.  An important property of these
products is that for any fixed
$a,b\in \goth A$ we have $a(n)b=0$ when $n$ is large enough. 
See  section \sec{defs} below for formal
definitions.

Conformal algebras where introduced in \cite{kac2} as the simplification
of the vertex algebra structure. Sometimes they are  also called 
``vertex Lie algebras''  \cite{dlm_poiss,primc}.  Curiously, similar algebraic
structures appeared in the Hamiltonian formalism in the theory of
non-linear evolution equations \cite{gd}. For more information on
conformal algebras see e.g. \cite{dk,kac2,kac_fd,zelmconf} and the references therein.

\subsection*{Formulation of the problem}
For any variety of algebras, like associative, Lie, Jordan, etc, 
there is the corresponding variety of conformal algebras, see
\sec{variety} for the rigorous statement. 
Given an associative conformal algebra
$\goth A$, we can define a different family of products 
$[n]:\goth A\times \goth A\to \goth A$ by 
\begin{equation}\label{fl:lieass}
a[n]b = a(n)b - \sum_{s\ge0}(-1)^{n+s}\frac1{s!}
D^s\big(b(n+s)a\big).
\end{equation}
With these new products (and the same derivation $D$),  the space $\goth A$ becomes a 
Lie conformal algebra, which we will denote by $\goth A^{(-)}$. This
is analogous to the fact that on any associative algebra $A$ one can
define the Lie algebra structure by using commutators, or Jordan
algebra structure by using anti-commutators. 
Like in these classical cases, a natural question is whether any 
Lie conformal algebra can be obtained as a subalgebra of $\goth
A^{(-)}$ for some associative conformal algebra $\goth A$. In this
case, $\goth A$ is called {\em an enveloping} algebra of $\goth L$. 
 
As shown in \cite{universal}, the answer to the above question is
negative for the following reason. Let  $\goth L$ be a Lie conformal
algebra, generated by a set  $\cal G$. Consider all words
$w = g_1(n_1)\cdots (n_{l-1})g_l\in\goth L$ (with arbitrary order of
parentheses)  for $g_i\in \cal G$ and $n_i\in
\Z_+$. Suppose that there  is an integer $S(l)=S_{\goth L,\cal G}(l)$ with the following
property: any word $w$ as above with $\sum_i n_i\ge S(l)$ is
zero. We call $S(l)$ {\em the locality
function} of $\goth L$. It is easy to show that if
$|\cal G|<\infty$, then $S_{\goth L,\cal G}(l)$ always exists and for a different set of
generators $\cal H$, the difference 
$|S_{\goth L,\cal G}(l)-S_{\goth L,\cal H}(l)|$
has at most linear growth in $l$. 
We have shown in \cite{universal}
that if a finitely generated Lie conformal algebra $\goth L$ is embeddable into an associative conformal algebra,
then $S(l)$ must have linear growth. 
On the other hand, for a free Lie conformal algebra the growth of
$S(l)$ is quadratic \cite{universal,cfva}, so it is  not embeddable
into associative. 

\subsection*{Results of this paper}
We prove the following theorem:
\begin{Thm}\label{thm:1}\sl
  Let $\goth L$ be a Lie conformal algebra, generated by a set $\cal
  G$. Assume that there is $K\in\Z$ such that  $S_{\goth L,\cal G}(l)\le K$ for
  any $l$. Then $\goth L$ is embeddable into  an associative
  conformal algebra $\goth A$ such that $S_{\goth A,\cal G}(l)\le K$
for all $l$. 
\end{Thm}

In particular, the condition of \thm{1} applies when the  algebra
$\goth L$ is nilpotent. We say that $\goth L$ is nilpotent of index
$k$ if all words $a_1(m_1)\cdots (m_{l-1})a_l$
are zero when $l\ge k$. Then we can prove
\begin{Thm}\label{thm:2}\sl
  Any nilpotent Lie conformal algebra has a nilpotent
  associative conformal enveloping algebra of the same index.
\end{Thm}

We remark that  not all algebras that satisfy the assumption of
\thm{1} are nilpotent. For example, the loop algebra $\goth L$ 
(see \sec{loop} below)
has locality function 1, but is not nilpotent. See \sec{bounded} for
other examples.

Both Theorems \ref{thm:1} and \ref{thm:2} will be derived from
\prop{main} in \sec{bounded}. The argument was partially
inspired by \cite{burde}.

\subsection*{Open questions}
First of all, it remains unclear whether a linear growth of the locality function is a
sufficient condition for the embedding of a Lie conformal algebra into
an associative. My guess is that this is false, but I don't have a
counterexample. 

An important special case is when the Lie conformal algebra 
is  of {\em finite type}, which means that it is a
module of finite rank over the  algebra of polynomials in $D$ \cite{dk}. It is
shown in  \cite{universal} that such algebras have linear locality functions.
It has been
conjectured \cite{universal} that a finite type Lie conformal algebra has a finite type
conformal associative enveloping algebra.  This conjecture is closely
related to another conjecture that states that  a finite type
torsion-free Lie conformal algebra
always has a faithful module of finite type  (see \sec{CEnd} for the
definitions and further discussion). 
This is the conformal analogue of classical Ado's theorem
\cite{burde,jacobson_lie}.

The results of
this paper show embeddability when the locality function is uniformly
bounded. In \sec{bounded} we conjecture that a central extension of 
a finite type Lie conformal algebra with bounded locality function
also has a bounded locality function. This would imply a general way
of generating examples of such algebras. 
\section{Definitions and Notations}\label{sec:defs}
All algebras and spaces are assumed to be over a field $\k$ of
characteristic zero. Throughout the paper we will use the divided
powers notation $x^{(n)} = \frac 1{n!}x^n$, and  $\Z_+$ will stand for
the set of non-negative integers.
 
\subsection{The definition}
\begin{Def}\cite{kac2}
A {\em conformal algebra} is a space $\goth A$ equipped by a
collection of bilinear products $A\otimes A \to A$,\  $a\otimes b
\mapsto a(n)b$, indexed by $n\in\Z_+$, and a linear map $D:\goth A \to
\goth A$ such that 
\begin{itemize}
\item[\hypertarget{C1}{(C1)}]
$a(n)b = 0$ for $n\gg 0$;
\item[\hypertarget{C2}{(C2)}]
$D\big(a(n)b\big) = (Da)(n)b + a(n)(Db) = -n\, a(n-1)b +  a(n)(Db)$. 
\end{itemize}
If only the condition \hyperlink{C2}{(C2)} is satisfied, then we will
call $\goth A$ a {\em preconformal algebra}.
\end{Def}

Iterating \hyperlink{C2}{(C2)}, we get 
\begin{equation}\label{fl:D}
\begin{gathered}
\big(D^{(k)}a\big)(n)b = (-1)^n \binom nk a(n-k)b,\\
a(n)\big(D^{(k)}b\big) =\sum_{s\ge 0} \binom ns D^{(k-s)}\big(a(n-s)b\big).
\end{gathered}
\end{equation}

\subsection{Formal series}
A typical way of constructing a conformal algebra is as follows. Take
a ``usual'' algebra $A$. Consider the space of formal series
$A[[z,z\inv]]$. We will write a series $\alpha \in A[[z,z\inv]]$ as 
$$
\alpha(z) = \sum_{n\in\Z}\alpha(n)\,z^{-n-1},\quad
\alpha(n)\in A.
$$

For $\alpha,\beta \in A[[z,z\inv]]$ define their
$n$-th product 
$(\alpha(n)\beta)(z) = \sum_{m\in\Z}\big(\alpha(n)\beta\big)(m)\,z^{-m-1}$ by
$$
(\alpha(n)\beta)(z) = \op{Res}_w \alpha(w)\beta(z)(z-w)^n,
$$
so that 
\begin{equation}\label{fl:serprod}
  \big(\alpha(n)\beta\big)(m) = \sum_{s=0}^n(-1)^s\binom ns \alpha(n-s)\beta(m+s).
\end{equation}

Series $\alpha,\beta \in A[[z,z\inv]]$ are {\em local} of
order $N\in \Z_+$, if 
$$
\alpha(w)\beta(z)(w-z)^N = 0.
$$
In terms of coefficients this means 
$$
\sum_{s=0}^N(-1)^s\binom Ns \alpha(n-s)\beta(m+s) = 0
$$
for any $m,n\in\Z$. 
It is easy to see that if $\alpha$ and $\beta$ are local of order $N$,
then $\alpha(n)\beta = 0$ for $n\ge N$. 

\begin{Prop}\label{prop:series}\cite{kac2,kac_fd}\sl\ 
 If $\goth A \subset A[[z,z\inv]]$ is a space of series such that 
 \begin{itemize}
 \item[(i)] 
any two series $\alpha,\beta \in \goth A$ are local;
\item[(ii)]  
$\alpha(n)\beta\in \goth A$ for any $\alpha,\beta \in \goth A$;
\item[(iii)]  
$\partial_z\alpha\in\goth A$ for any  $\alpha\in \goth A$,
 \end{itemize}
then $\goth A$ is a conformal algebra with $D = \partial_z$.
If only conditions {\rm (ii)} and {\rm (iii)} hold, then $\goth A$ is
a preconformal algebra. 
\end{Prop}
\begin{Rem}
If $A$ is either associative or Lie algebra, then condition (i) of
\prop{series} can be weakened: it is enough to assume the locality
only for a set of generators $\cal G$ of $\goth A$. This fact is know
as  {\em the Dong's lemma}.
\end{Rem}

\subsection{Coefficient algebra}
Conversely, any conformal algebra can be obtained as in
 \prop{series}. 
Moreover, to any conformal algebra $\goth A$ there corresponds a ``usual''
algebra $A=\cff \goth A$, called the {\em coefficient algebra} of
$\goth A$, and the inclusion $\pi:\goth A \hookrightarrow A[[z,z\inv]]$
with the following universal property. For any other 
homomorphism $\varphi:\goth A \to B[[z,z\inv]]$ of $\goth A$ to the
space of formal series, such that $\varphi(\goth A)$ satisfies the
conditions of \prop{series}, there is a unique algebra homomorphism
$\rho:A\to B$ such that $\rho(\pi(a)) = \varphi(a)$ for any $a\in
\goth A$. 

The coefficient algebra  $A = \cff\goth A$ is constructed in the
following way. Consider the space of Laurent series 
$\goth A[t,t\inv]$ in an independent variable $t$ with coefficients in
$\goth A$. For $a\in \goth A$,
denote $a(n) = at^n$.  As a linear space $A$ is isomorphic to the
quotient of $\goth A[t,t\inv]$ over the subspace generated by 
the vectors $(Da)+n\,a(n-1)$ for $a\in \goth A$. The formula for the product in $A$ is
derived from \fl{serprod}:
$$
a(m)b(n)=\sum_{s\ge 0} \binom ms \big(a(s)b\big)(m+n-s).
$$
Note that the sum here is finite due to \hyperlink{C1}{(C1)}.
The canonical inclusion $\goth A \to A[[z,z\inv]]$ is given by 
$a\mapsto \sum_{n\in\Z}a(n)\,z^{-n-1}$.

In \sec{proof} below we are going to need to following fact:
\begin{Lem}\label{lem:coeff}\cite{kac_fd,freecv}\ \sl
Assume that $\goth A$ is a free $\k[D]$-module, and let $\cal B\subset
\goth A$ be its basis over $\k[D]$. Then the set
$\bigset{b(n)}{b\in\cal B,\,n\in\Z}$ is a $\k$-linear basis of
$\cff\goth A$. 
\end{Lem}

\begin{Rem}
The requirement that $\goth A$ is a free $\k[D]$-module is not as
restrictive as it might appear. In any conformal algebra $\goth A$
one can define the so-called {\em torsion ideal}  
$\goth t = D\!\op{-tor}\goth A + \bigcap_{n\ge 0} D^n \goth A$, where 
$D\!\op{-tor}\goth A = \set{a\in\goth A\,}{\,\exists \, p(D)\in\k[D], \,
  p(D)\neq 0,\,
p(D)a=0}$ is the $D$-torsion of $\goth A$, so that $\goth A/\goth t$
is a free $\k[D]$-module. It is easy to show \cite{dk} that $\goth t$
belongs to the left annihilator of $\goth A$, i.e. $a(n)b=0$ for any 
$a\in\goth t$, \ $b\in\goth A$ and $n\in \Z_+$.   
\end{Rem}

\subsection{Varieties of conformal algebras}\label{sec:variety}
In the case when  the coefficient algebra
$A=\cff\goth A$ belongs to a certain variety of algebras, the conformal algebra
$\goth A$ is said to belong to the corresponding conformal variety. 
For example, if $A$ is
an associative (respectively, a Lie or a Jordan) algebra, then 
$\goth A$ is called an associative conformal (respectively, a  Lie
conformal or a Jordan conformal) algebra. Moreover, if $A$ belongs to
a certain variety of algebras, then any conformal subalgebra of
$A[[z,z\inv]]$ as in \prop{series} belongs to the corresponding
conformal variety. In this paper we deal only
with Lie or associative conformal algebras. To distinguish between
them, we will denote the products in a  conformal algebra by
$[n]$ whenever the product in the coefficient algebra is denoted by
the brackets $\ad{\,\cdot}{\cdot}$. 
 
There is a correspondence between the identities in a conformal algebra
and the identities in its coefficient algebra. An identity $R$ holds
in $A =\cff\goth A$ if and only if a certain identity (or family of
identities) $\op{Conf}R$ holds in $\goth A$. For example, the
associativity $(ab)c=a(bc)$, the Jacoby identity $\ad{\ad ab}{c} =
\ad{a}{\ad bc} - \ad{b}{\ad ac}$ and the (skew-)symmetry
$ab=\pm \,ba$ correspond to the following conformal identities respectively:

conformal associativity:
\begin{equation}\label{fl:assoc}
  \big(a(m)b\big)(n)c = \sum_{s\ge 0} (-1)^s \binom ms
  a(m-s)\big(b(n+s)c\big)
\end{equation}

conformal Jacoby identity:
\begin{equation*}
  \big(a[m]b\big)[n]c = \sum_{s\ge 0} (-1)^s \binom ms
  \Big(a[m-s]\big(b[n+s]c\big) - b[n+s]\big(a[m-s]c\big)\Big)
\end{equation*}

quasi-symmetry:
\begin{equation}\label{fl:qs}
  a(n)b  = \pm\sum_{s\ge 0} (-1)^{s+n}D^{(s)}\big(b(n+s)a\big)
\end{equation}

We will need in \sec{proof} the following strengthening of the above correspondence.
\begin{Lem}\label{lem:preconf}\sl
Let $A$ be an associative algebra, and  $\goth A\subset A[[z,z\inv]]$
a preconformal algebra of formal series with coefficients in $A$. Then
the identity \fl{assoc} holds in $\goth A$. 
\end{Lem}
\begin{proof}
Take $k\in\Z$. The $k$-th coefficient of the left- and right-hand sides
of \fl{assoc} are, respectively  
\begin{equation*}
\sum_{i,j\ge 0} (-1)^{i+j} \binom mi\binom nj a(m-i)b(n+i-j)c(k+j) 
\end{equation*}
and 
\begin{equation*}
\sum_{i,j,s\ge 0} (-1)^{i+j+s} \binom ms\binom {m-s}i\binom {n+s}j 
a(m-s-i)b(n+s-j)c(k+i+j).
\end{equation*}
Replace the indices in the second formula by the rule $i\to i-s$, \
$j\to j-i+s$, 
and then we are done by the combinatorial identity
$$
\sum_{s=0}^i (-1)^{s+j}\binom ms \binom {m-s}{i-s}\binom {n+s}{j-i+s} 
=  (-1)^{i+j} \binom mi\binom nj.
$$
\end{proof}

For a Lie conformal algebra $\goth L$ denote by $Z(\goth L) =
\set{a\in\goth L}{a(n)b=0\ \forall\, b\in\goth L,\,n\in\Z_+}$ the
center of $\goth L$. Due to the quasi-symmetry \fl{qs}, we have 
$b(n)a = 0$ for any $a\in Z(\goth L)$, \ $b\in\goth L$ and $n\in\Z_+$.

\subsection{The relation between associative and Lie conformal
  algebras} \label{sec:lieass}
Let $\goth A$ be an associative conformal algebra. Then we can define
another family of products $[n]$, \ $n\in\Z_+$, on $\goth A$ by the 
formula \fl{lieass}. This will define a Lie conformal algebra
structure on $\goth A$, which we will denote by $\goth
A^{(-)}$. Recall that any associative algebra $A$ can be turned into a
Lie algebra $A^{(-)}$ by taking the commutator $\ad ab = ab-ba$ for
the product. It is easy to check that $(\cff \goth A)^{(-)} = \cff (\goth A^{(-)})$.

Here is another useful formula that holds in $\goth A$:

\begin{equation}\label{fl:adconf}
  a(m)b(n)c-b(n)a(m)c = \sum_{s\ge 0} \binom ms \big(a[s]b\big)(m+n-s)c
\end{equation}

In \sec{proof} we will deal with the following situation. Let $\goth
L$ be a Lie conformal algebra and $L= \cff\goth L$ be its coefficient 
Lie algebra. Let $A\supset L$ be an associative enveloping algebra of
$L$. Then we get $\goth L \subset A[[z,z\inv]]$. Let $\goth A \subset
A[[z,z\inv]]$ be the associative preconformal algebra generated by
$\goth L$. By \lem{preconf}, the conformal associativity \fl{assoc}
holds in $\goth A$. The following statement is checked in a similar
way.

\begin{Lem}\label{lem:adconf}\sl
The formula \fl{adconf} holds for any $a,b\in\goth L$, \  $c\in\goth A$
and $m,n\in\Z_+$.
\end{Lem}

\subsection{Example: loop algebras}\label{sec:loop}
Let $\goth g$ be an  algebra. Let $L = \goth g[[t,t\inv]$. Denote
$a(m) = at^m$ for $a\in\goth g$ and $m\in\Z$, so that $a(m)b(n)
  = (ab)(m+n)$. For $a\in\goth g$ set 
$$
\~a =
  \sum_{m\in\Z}a(m)\,z^{-m-1}\in L[[z,z\inv]].
$$ 
It is easy to see that 
$\~a$ and $\~b$ are local of order 1 and $\~a(0)\~b = \~{ab}$.
Let $\goth L \subset L[[z,z\inv]]$ be the conformal algebra
generated by $\~a$ for all $a\in\goth g$. It is called the {\em loop
  algebra} of $\goth g$.  As a $\k[D]$-module, 
$\goth L$ is freely generated by $\cal G = \set{\~a}{a\in\goth g}$, so 
\lem{coeff} implies that  $L=\cff\goth L$.

We remark that a monomial $\~a_1(m_1)\cdots
\~a_{l-1}(m_{l-1})\~a_l\in \goth L$ (with  arbitrary order of 
parentheses) is equal to 0 if $\sum m_i >0$, so the locality function
$S_{\goth L,\cal G}=1$. It is easy to show that if a conformal algebra
has locality function 1, then it must be a loop algebra.

In the case when $\goth g$ is a Lie algebra, it often comes with an
invariant bilinear form $\formdd$, and then the corresponding loop
algebras $L$ and $\goth L$ have central extensions $\^L = L\oplus
\k\c$ and $\^{\goth L} = \goth L \oplus \k\c$. The brackets in $\^L$
are  $\ad{a(m)}{b(n)} 
=\ad ab(m+n) + m \delta_{m,-n} \form ab$, the locality of $\~a$
and $\~b$ is 2, and the conformal products are 
$\~a(0)\~b = \~{\ad ab}$, \ $\~a(1)\~b =\c$. We identify $\c$ with
$\c(-1)\in\^L[[z,z\inv]]$. This is called the {\em affine Lie
  algebra} corresponding to $\goth g$. The locality function of 
$\^{\goth L}$ corresponding to the generators $\cal G \cup \{\c\}$ 
is equal to 2. 

\subsection{Representations of conformal algebras}\label{sec:CEnd}
Let $M$ be a $\k[D]$ module of finite rank. 
\begin{Def}\cite{ck}  
A {\em conformal operator} $\alpha$ on $M$ is a series 
$\alpha = \sum_{n\ge 0} \alpha(n)\,z^{-n-1} \in \op{End}(M)[[z\inv]]$,
such that 
\begin{itemize}
\item[(CO1)] for any fixed $v\in M$ we have $\alpha(n)v = 0$ for $n\gg
  0$,
\item[(CO2)] $\ad D{\alpha(n)} = -n\,\alpha(n-1)$. 
\end{itemize}
Denote the space of all conformal operators by $\op{CEnd}(M)\subset
\op{End}(M)[[z\inv]]$. 
\end{Def}

In fact, any conformal operator $\alpha \in \op{CEnd}(M)$ is 0 on 
$D\!\op{-tor} M$, so we can assume that $M$ is a free $\k[D]$-module
without loss of generality. 

We observe that the formula \fl{serprod} makes sense when 
$\alpha,\beta \in \op{CEnd}(M)$, and also $\op{CEnd}(M)$
is closed under the derivation $D = \partial_z$, so 
it can be shown that $\op{CEnd}(M)$ is an associative conformal
algebra \cite{ck}. 

Let $\goth A$ be an associative (respectively, a Lie) conformal algebra,
then by definition, $M$ is a module over $\goth A$ if there is a
conformal algebra homomorphism $\goth A \to \op{CEnd}(M)$
(respectively, $\goth A \to \op{CEnd}(M)^{(-)}$). For example, the
algebra $\goth A$ is a module over itself with the representation map 
$\goth A \to \op{CEnd}(\goth A)$ given
by  $a\mapsto \sum_{n\ge 0} a(n)\,z^{-n-1}$.  

It follows that if a Lie conformal algebra $\goth L$ has a faithful
finite type module $M$, then $\op{CEnd}(M)$ is an associative
conformal enveloping algebra of $\goth L$. In particular, this applies
to the case when $Z(\goth L)=0$ so that $\goth L$ is a faithful module
over itself. Together with Theorems \ref{thm:1} and \ref{thm:2} this
provides some grounds to the conjecture that any finite type Lie conformal
algebra is embeddable into an associative conformal algebra. 

For further information about these and other conformal algebras
consult e.g. the reviews \cite{kac_fd,zelmconf} and the references therein.
\section{Conformal algebras with bounded locality  function}
Let $\goth L$ be a conformal algebra generated by
a set $\cal G\subset\goth L$. Recall that the {\em locality function}
$S(l) =S_{\goth L,\cal G}(l)$ is an integer such that any word
$w = g_1(n_1)\cdots (n_{l-1})g_l\in\goth L$, where 
$g_i\in \cal G$ and $n_i\in
\Z_+$, is zero whenever $\sum_i n_i \ge S(l)$. If $|\cal G|<\infty$,
then the existence of $S(l)$ follows from 
\hyperlink{C1}{(C1)}. If $\goth L$ is a Lie conformal algebra, then 
the quantitative version of Dong's lemma \cite{cfva} implies 
that if the locality $N(a,b)$ of any two generators 
$a,b\in\cal G$ is uniformly bounded by $N\in\Z_+$, then 
$S(l) \le \frac 12 N l(l-1) -l +1$.

\subsection{The main Proposition}\label{sec:bounded}
To every generator $g\in \cal G\subset \goth L$ we assign a weight $\wt g
\in \Z_+$. (The word ``degree'' will be used later for different
purpose). 
For a monomial $w=g_1[n_1]\cdots [n_{l-1}]g_l$ or 
$w=g_1(n_1)\cdots (n_{l-1})g_l$ (with
arbitrary order of parentheses), where $g_i\in \cal G$ and $n_i\in
\Z_+$, we set $\wt w = \sum_i \wt g_i + \sum_i n_i$. We also
set $\wt D = -1$.

Both \thm{1} and\thm{2} follow from the following statement:
\begin{Prop}\label{prop:main}\sl
Assume that there is an integer $r\ge 0$ such that any monomial 
of weight $r$ or more is equal to zero in $\goth L$. Then 
$\goth L$ is embedded into an enveloping conformal associative algebra
$\goth A$, that has the same property: any conformal monomial $w$ in 
$\cal G$ is zero in $\goth A$ whenever $\wt w \ge r$. 
\end{Prop}

\thm{1} (respectively, \thm{2}) 
is a special case of \prop{main} obtained by setting $\wt g =
1$ (respectively, $\wt g = 0$) for every $g\in \cal G$.  We will prove
\prop{main} in \sec{proof}.

We have seen one example of a Lie conformal algebra with bounded
locality function in \sec{loop}. It is easy to see that any central
extension of a loop algebra also has a bounded locality function. We
state the following conjecture.

\begin{Conj}\sl
  Let $\goth L$ be a Lie conformal algebra of finite type generated by
  a finite set $\cal G$ so that $S_{\goth L, \cal G}(l)<K$. Let $M$ be
  a $\k[D]$-module, generated by a set $\cal M$, on which $\goth L$
  acts trivially, and let $\^{\goth L}$
  be a central extension of $\goth L$. Then 
$S_{\^{\goth L},\,\cal G\cup \cal M}$ is also uniformly bounded.
\end{Conj}

By ``central extension'' we mean that there is a short exact sequence
of conformal algebra homomorphisms
$0 \to M \to \^{\goth L}\to \goth L \to 0$, 
such that $M\subseteq Z(\^{\goth L})$.

\subsection{Proof of \prop{main}}\label{sec:proof}
\subsubsection{Two filtrations on $\goth L$}\label{sec:gothL}
Let $\goth L$ be a Lie conformal algebra, satisfying the conditions of\break
\prop{main}.  Define  a filtration 
$$
\goth L\supseteq \ldots\supseteq \goth L'_{i-1}\supseteq \goth L'_i
\supseteq \goth L'_{i+1} \supseteq\ldots \supseteq \goth L'_r = 0
$$
on $\goth L$ by setting 
$$
\goth L'_i = \spn\bigset{w = D^m g_1[n_1]\cdots [n_{l-1}]g_l}{
g_j\in \cal G,\ \wt w  \ge i}.
$$
We have $\goth L'_i = 0$ for $i\ge r$ due to the fact that there are no words
of weight $l$ or more in $\goth L$.   Clearly, we also 
 have $\bigcup_i \goth L'_i = \goth L$.
For an element $a\in\goth L$ set $\deg' a = \max\set i{a\in \goth
  L'_i}$. Note that for a Lie conformal monomial $w$ in $\cal G$ 
we have $\wt w \le \deg' g$.  

Here are  some easy properties of this filtration that we are
going to need:

\begin{Lem}\label{lem:Lpr}\sl
  \begin{enumerate}
  \item\label{Lpr:ideal} 
$\goth L'_i [n]\goth L'_j \subseteq \goth L'_{i+j+n}$, \ 
$D\goth L'_i \subseteq \goth L'_{i-1}$.
  \item\label{Lpr:L0} 
$\goth L = \k[D]\goth L'_0$.
\end{enumerate}
\end{Lem}
\begin{proof}
(\ref{Lpr:ideal}) follows from the fact that the formulas \fl{D}
 are homogeneous. 
To prove (\ref{Lpr:L0}), note that since every
  generator $g\in\cal G$ has $\wt g \ge 0$, an element $a\in\goth L$
  of negative degree must belong to $D\goth L$.    
\end{proof}

Next we set 
$$
\goth L_i = \bigset{a\in\goth L}{\exists n: D^n a\in\goth
  L'_{i-n}}.
$$ 
This defines another filtration on $\goth L$ of the form
$$
\goth L\supseteq \ldots\supseteq \goth L_{i-1}\supseteq \goth L_i
\supseteq \goth L_{i+1} \supseteq\ldots .
$$
Set $\goth L_\infty = \bigcap_i \goth L_i$ and $\deg a =
\sup\set{i}{a\in\goth L_i}$ for $a\in\goth L$. Clearly, we have $\goth
L_i'\subseteq \goth L_i$, therefore $\deg a \ge \deg' a$. More
precisely, we have 
\begin{equation}\label{fl:deg}
 \deg a = \sup_n\{n+\deg'D^na\}. 
\end{equation}

 Here are some properties of the filtration $\{\goth L_i\}$:

\begin{Lem}\label{lem:Lconf}\sl
  \begin{enumerate}
  \item\label{Lconf:ideal} 
$\goth L_i[n]\goth L_j \subseteq \goth L'_{i+j+n}$. In particular,
$\goth L_i[n]\goth L_j = 0$ if $i+j+n\ge r$.
  \item\label{Lconf:Z}
 $\goth L_r\subseteq Z(\goth L)$.
  \item\label{Lconf:D}
$\deg Da = \deg a - 1$ for any $a\in\goth L$.
  \end{enumerate}
\end{Lem}

\begin{proof} (\ref{Lconf:ideal})\ 
Let $a\in\goth L_i$ and $b\in\goth L_j$. Since $a[n]b = 0$ for $n\gg
0$, we can, using induction,  assume that $a[s]b\in\goth L'_{i+j+s}$
for any $s>n$. 

There are $k,m\in\Z_+$, such that $D^{(k)}a\in\goth L'_{i-k}$ and
 $D^{(m)}b\in\goth L'_{j-m}$. Using \lemm{Lpr}{ideal} and \fl{D}, we get
\begin{align*}
  \goth L'_{i+j+n}\ni \big(D^{(k)}a\big)[k+m+n]\big(D^{(m)}b\big) &=
(-1)^k \binom{k+m+n}{k}\sum_{s=0}^m \binom{m+n}{m-s}\, D^{(s)}\big( a[n+s]b\big)\\
&\equiv  (-1)^k \binom{k+m+n}k\binom{m+n}{m} \ a[n]b \mod \goth L'_{i+j+n},
\end{align*}
since by \lemm{Lpr}{ideal} and induction, 
$D^{(s)}\big( a[n+s]b\big) \in \goth L'_{i+j+n}$ for $s>0$.

\smallskip\noindent(\ref{Lconf:Z})\ 
Since $\goth L'_0\subseteq \goth L_0$, \lemm{Lpr}{L0} implies that 
$\goth L = \k[D]\goth L_0$, and by  (\ref{Lconf:ideal}) we have $a[n]b = 0$
for any $a\in\goth L_0$ and $b\in\goth L_r$. 

\smallskip\noindent(\ref{Lconf:D})\ 
By \lemm{Lpr}{ideal} we have $\deg'Da \ge \deg'a-1$. Using this and \fl{deg}, we get
$$
\deg Da = \sup_{n\ge 0}\big\{n+\deg'D^{n+1}a\big\} = 
\sup_{n\ge 0}\big\{n-1+\deg'D^na\big\} = \deg a-1.
$$
\end{proof}
It follows from (\ref{Lconf:Z}) and (\ref{Lconf:D}) that 
$\goth L_\infty$ is a central ideal of $\goth L$.

\subsubsection{The basis $\cal B$}\label{sec:B}
Let $\cal B_i \subset \goth L_i$ be a 
$\k$-linear basis of $\goth L_i$
modulo $\goth L_{i+1} + D\goth L_{i+1}$. By \lemm{Lpr}{L0}, if $i<0$, then 
$\cal B_i = \varnothing$. Denote $\cal B = \bigcup_{i=0}^{r-1}
\cal B_i$. Let $\goth T = \k[D]\goth L_r$. By \lemm{Lconf}{ideal} and (\ref{Lconf:Z}),
this is  a central ideal of $\goth L$. 

\begin{Lem}\label{lem:B}\sl
The set $\cal B$ is a $\k[D]$-linear basis of $\goth L \mod \goth T$. 
The expansion of an element $a\in \goth L$ of $\deg a = i$
in this basis is
\begin{equation}\label{fl:aexp}
a=\sum_{n\in\Z_+,\,b\in\cal B}k_{n,b}\,D^nb+ a_0,\qquad
k_{n,b}\in\k, \ a_0\in \goth T,
\end{equation}
such that $\deg D^nb = \deg b - n \ge i$ whenever $k_{n,b}\neq 0$.
\end{Lem}
\begin{proof}
First we show that  any element $a\in \goth L_i$ has expansion
\fl{aexp}. Indeed, if $i\ge r$, this is obvious, so by induction we
can assume that any $a\in\goth L_{i+1}$ has such an expansion. Now, 
we can decompose $a=a_1+a_2$ so that $a_1 \in \spn_\k\cal B_i$ and 
$a_2\in \goth L_{i+1} + D\goth L_{i+1}$; by induction, $a_2$ has an
expansion \fl{aexp}, therefore, so does $a$.

This shows that $\cal B$ spans $\goth L$ modulo $\goth
T$ over $\k[D]$. Let us prove that the set $\cal B$ is linearly
independent over  $\k[D]$ modulo $\goth T$. 

We observe that in the $\k[D]$-module $\goth L/\goth T$ we have $\ker
D = 0$. Indeed, every element $t\in\goth T$ can be written as $t =
Dt_0+t_1$ for $t_0\in\goth T$ and $t_1\in\goth L_r$. So if $Da
\in\goth T$, then write $Da = Dt_0+t_1$, and get $D(a-t_0) =
t_1\in\goth L_r$, hence by \lemm{Lconf}{D} we get $\deg (a-t_0) \ge
r+1$.  Therefore, $a-t_0\in\goth T$, hence also $a\in\goth T$.

Now assume that we have a linear relation
$$
\sum_{n\in\Z_+,\,b\in\cal B} k_{n,b}\, D^nb\in\goth T,\qquad
k_{n,b}\in\k.
$$
By the previous paragraph, we can assume that the set 
$\cal B'=\set{b\in\cal B}{k_{0,b}\neq 0}$ in non-empty.  
Let $i$ be the minimal degree of the elements in this set, and
let $\cal B'_i=\cal B_i\cap \cal B'$ be the elements degree $i$ in $\cal B'$. 
The linear relation above implies that $\cal B'$ is linear dependent modulo
$D\goth L + \goth L_{i+1}$, which contradicts the definition of $\cal
B_i$, since $D\goth L \cap \goth
L_i = D\goth L_{i+1}$  due to \lemm{Lconf}{D}. 
\end{proof}

\subsubsection{The algebra $\goth U$}\label{sec:gothU}
Let $L = \cff\goth L$ be the coefficient Lie algebra of $\goth
L$. Consider its universal enveloping algebra $U(L)$ and let $U$ be
the augmentation ideal of $U(L)$.

Consider the space  $U[[z,z\inv]]$ of formal series with coefficients
in $U$.  Let 
$\goth U \subset U[[z,z\inv]]$ be the associative preconformal algebra
generated by the series $\sum_{m\in \Z}a(m)\,z^{-m-1}$ for $a\in\goth
L$. By \lem{preconf} and \lem{adconf} the identity  \fl{assoc} holds
for any $a,b,c\in\goth U$ and 
\fl{adconf} holds  for any $a,b\in\goth L$ and  $c\in\goth U$.

Define a filtration 
$
\goth U\supseteq \ldots\supseteq \goth U_{i-1}\supseteq \goth U_i
\supseteq \goth U_{i+1} \supseteq\ldots 
$
 on $\goth U$ by setting
$$
\goth U_i = \spn\Bigset{D^m a_1(m_1)\cdots a_{l-1}(m_{l-1})a_l}{
a_j\in \goth L,\ m_j\in\Z_+,\ \sum_{j=1}^{l-1} m_j + \sum_{j=1}^l \deg a_j -m \ge i}.
$$
Here the order of parentheses is arbitrary, but note that using the
formula \fl{assoc}, it is enough to take only right-normed words.  

The filtration $\{\goth U_i\}$ is defined in a similar way to the
filtration $\{\goth L'_i\}$ on $\goth L$, constructed in \sec{gothL},
so it satisfies the properties, analogous to \lemm{Lpr}{ideal}, which are  
proved in the same way:

\begin{Lem}\label{lem:Uconf}\sl
$\goth U_i (n)\goth U_j \subseteq \goth U_{i+j+n}$, \ 
$D\goth U_i\subseteq \goth U_{i-1}$.
\end{Lem}

\subsubsection{}\label{sec:T}
Let $T = \cff\goth T \subset L$ be the coefficient algebra of $\goth
T$. Then $T$ is a central ideal of $L$ and we have 
$L/T = \cff (\goth L/\goth T)$. 
By \lem{coeff}, since $\cal B$ is a $\k[D]$-linear basis of $\goth L \mod
\goth T$, the set $\set{b(n)}{b\in\cal B,\, n\in\Z}\subset L$ is a $\k$-linear
basis of $L \mod T$.

Let $N\subset U$ be the ideal of $U$ generated by $T = \cff \goth T$.
Then the algebra $U/N$ is equal to the augmentation ideal of the
universal enveloping algebra $U(L/T)$.  Choose a linear order
on $\cal B$. For $m,n\in\Z$ and $a,b\in \cal B$ we will write 
$a(m)<b(n)$ if either $m<n$ or $m=n$ and $a<b$. The PBW theorem
states that the set 
\begin{equation}\label{fl:S}
\cal S=\bigset{b_1(n_1)\cdots b_l(n_l)}{b_i\in\cal B,
\, n_i\in\Z, \, b_i(n_i)\ge b_{i+1}(n_{i+1})}\subset U
\end{equation}
is a $\k$-linear basis of $U \mod N$. We will order the words from
$\cal S$ first by length and then alphabetically from right to left,
so that for  
$u=b_1(n_1)\cdots b_l(n_l),\, u'=b'_1(n'_1)\cdots b'_l(n'_l)\in \cal S$ 
we write $u\le u'$ if
$b_l(n_l)=b'_l(n'_l), \ldots, b_{i+1}(n_{i+1})=b'_{i+1}(n'_{i+1})$,
but $b_i(n_i)\le b'_i(n'_i)$ for some $1\le i\le l$.

Since
$T$ is a central ideal of $L$, the set $N_2 = TN \subset N$ is a
proper subideal of $N$ spanned by all words $a_1(m_1)\cdots a_l(m_l)$
for $a_i\in \goth L$, \ $m_i\in \Z$, \  $l\ge 2$, such that
$a_i\in\goth T$ for some $1\le i\le l$. Clearly we have $N_2\cap L =
0$.

Similarly, we define $\goth N\subset \goth U$ to be the span of all words 
$D^n a_1(n_1)\cdots a_{l-1}(n_{l-1})a_l$, \ $a_i\in\goth L$ 
(for arbitrary order of parentheses), such that $a_i\in\goth T$ 
for some $1\le i\le l$, and 
$\goth N_2\subset \goth N$ to be the span of the same words of length at
least 2. 
Clearly, $\goth N$ and $\goth N_2$ are ideals of $\goth U$ such that 
$u(n) \in N$ for $u\in \goth N$ and $u(n)\in N_2$ for $u\in\goth N_2$,
and we have   $\goth N_2 \cap\goth L = 0$.

\subsubsection{Basis in $\goth U$}
Recall from \sec{B} that we have a set $\cal B\subset \goth L$ which
is a $\k[D]$-linear basis of $\goth L/\goth T$. Define
$$
\cal W = \left\{b_1(n_1)\big(b_1(n_2)\cdots\big(
  b_{l-1}(n_{l-1})b_l\big)\hskip-2pt\cdot\hskip-2pt\cdot\hskip-2pt\cdot\hskip-3pt\big)\in\goth U\,\left|\,
\begin{gathered}
b_i\in\cal B, \ n_i\in\Z_+\\
  b_i(n_i)\ge b_{i+1}(n_{i+1})\ \text{for}\ 1\le i\le l-2
\end{gathered}
\right.\right\}.
$$
For  $w = b_1(n_1)\big(b_1(n_2)\cdots\big(b_{l-1}(n_{l-1})b_l
\big)\hskip-2pt\cdot\hskip-2pt\cdot\hskip-2pt\cdot\hskip-3pt\big)
\in\cal W$ set $\deg w = \sum_j \deg b_j + \sum_j n_j$, so that 
$w\in \goth U_{\deg w}$.

\begin{Lem}\label{lem:Ubasis}\sl
  The set $\cal W$ is a $\k[D]$-linear basis of $\goth U/\goth N$,
  such that the expansion of an element $u\in\goth U_i$ in this basis has
  form
  \begin{equation}\label{fl:uexp}
u = \sum_{n\in\Z_+,\,w\in\cal W} k_{n,w}\,D^{(n)}w+u_0,\quad
k_{n,w}\in \k,\ u_0\in \goth N,
  \end{equation}
where $\deg D^{(n)}w = \deg w - n \ge i$ whenever $k_{n,w}\neq 0$.
\end{Lem}

\begin{proof}
Let us show first that  $\cal W$ is linearly independent over $\k[D]$ 
modulo $\goth N$. Suppose 
$$
u = \sum_{n,w} k_{n,w}\,D^{(n)}w\in
\goth N.
$$ 
Then  $u(-1)= \sum_{n,w} k_{n,w}\,w(-n-1) \in N$.  For 
$w = b_1(n_1)\big(b_1(n_2)\cdots\big(
  b_{l-1}(n_{l-1})b_l\big)\hskip-2pt\cdot\hskip-2pt\cdot\hskip-2pt\cdot\hskip-3pt\big)
\in \cal W$ we compute, iterating \fl{assoc},
\begin{multline*}
 w(-n-1) = \sum_{i_1=0}^{n_1}\cdots \sum_{i_{l-1}=0}^{n_{l-1}}
(-1)^{i_1+\ldots+i_{l-1}} \binom{n_1}{i_1}\cdots \binom{n_{l-1}}{i_{l-1}}\\
\times b_1(n_1-i_1)\cdots b_{l-1}(n_{l-1}-i_{l-1})b_l(i_1+\ldots +i_{l-1}-n-1).
\end{multline*}
If we expand $w(-n-1)$ into a linear combination of the elements of
$\cal S$ (see \fl{S}) modulo
$N$, then the minimal term among the terms of maximal length in 
this expansion will be 
$$
b_1(n_1)\cdots
b_{l-1}(n_{l-1})b_l(-n-1).
$$
 We observe that these terms are different
for every pair $n\in \Z_+$,\, $w\in \cal W$. Therefore, the set 
$\set{w(-n-1)}{n\in \Z_+,\, w\in \cal W}$ is linearly independent modulo
$N$, hence all $k_{n,w}=0$.

Now we show that any element $u\in \goth U_i$ has expansion \fl{uexp}. 
By definition, $\goth U_i$ is spanned by words 
$D^m a_1(m_1)\cdots a_{l-1}(m_{l-1})a_l$ for $a_j\in\goth L$ and
$m_j\in\Z_+$, such that $\sum_j m_j + \sum_j\deg a_j -m\ge i$.
Expand every $a_j$ in such a word into a linear combination \fl{aexp},
and then use \fl{assoc} and \fl{adconf} to write this word in the form
\fl{uexp}. By \lem{B}, every term in the expansion \fl{aexp} for $a_j$
will have degree at least $\deg a_j$, so the condition on degrees in
\fl{uexp} follows from
the fact that the relations  \fl{assoc} and \fl{adconf} are homogeneous.
\end{proof}

\subsubsection{The ideal $\goth I$}
Set  
$$
\goth I = \spn_{\k[D]} \left\{a_1(m_1)\cdots a_{l-1}(m_{l-1})a_l\in
\goth U\ \left|
\ \begin{gathered}
  a_j\in\goth L,\  m_j\in\Z_+,\  l\ge 2\\
\textstyle{ \sum_j m_j +\sum_j \deg a_j \ge r}
\end{gathered}
\right.\right\}.
$$

\begin{Lem}\label{lem:I}\sl
  \begin{enumerate}
  \item\label{I:ideal} 
$\goth I$ is an ideal of $\goth U$.
  \item\label{I:loc} 
For any $a,b\in\goth L$ we have $a(n)b\in \goth I$ for $n\gg 0$.
  \item\label{I:cap} 
$\goth I\cap \goth L = 0$.
  \end{enumerate}
\end{Lem}
\begin{proof}
(\ref{I:ideal})\ 
Let $u = a_1(m_1)\cdots a_{l-1}(m_{l-1})a_l$ be a generator of $\goth I$. 
Then $a(n)u \in \goth I$ for any 
$a\in \goth L_0$ and $n\in\Z_+$. Since $\goth L = \k[D]\goth L_0$ by \lemm{Lpr}{L0}, 
we get  $\goth L(n)\goth I \subseteq \goth I$ for any 
$n\in\Z_+$. But $\goth U$ is generated by $\goth L$ as an algebra,
therefore \fl{assoc} implies that $\goth U(n)\goth I \subseteq \goth I$ for any 
$n\in\Z_+$.

\smallskip\noindent(\ref{I:loc})\ 
We have $a(n)b\in \goth I$ for $n \ge l-\deg a-\deg b$.

\smallskip\noindent(\ref{I:cap})\ 
Let $u = a_1(m_1) \cdots a_{l-1}(m_{l-1})a_l\in\goth I$ as above. 
It follows from \lem{Ubasis} and \sec{T} 
that the expansion \fl{uexp} of $u$
will have the following two properties:
\begin{itemize}
\item[(i)]
$\deg w \ge r$ whenever $k_{n,w}\neq 0$;
\item[(ii)]
$u_0\in \goth N_2$. 
\end{itemize}
Clearly, the expansion \fl{uexp} of any $\k[D]$-linear combination of such
elements $u$ also has properties (i) and (ii). We are left to
note, that if a linear combination \fl{uexp} with  properties
(i) and (ii) belongs to $\goth L$, then it is 0. Indeed,
any word $w\in\cal W$ of degree $r$ or more has length at least 2,
since $\deg b \le r-1$ for any $b\in\cal B$. Therefore, all $k_{n,b}
= 0$ and the combination is in $\goth N_2$. 
But we also have $\goth N_2 \cap \goth L = 0$.
\end{proof}

Now \prop{main} easily follows from \lem{I}: Take $\goth A = \goth U/\goth I$. This is
an associative conformal algebra, since the conformal associativity
holds by \lem{preconf} and the locality is due to (\ref{I:loc}), it contains
$\goth L$ because of (\ref{I:cap}) and any word $w=g_1(n_1)\cdots
(n_{l-1})g_l$, \ $g_i\in\cal G$, of weight $r$ or more belongs to
$\goth I$, therefore, $w=0$ in $\goth A$.

\bibliography{../vertex,../my,../general,../conformal}

\def\cprime{$'$}
\begin{thebibliography}{10}

\bibitem{burde}
D.~Burde.
\newblock On a refinement of {A}do's theorem.
\newblock {\em Arch. Math. (Basel)}, 70(2):118--127, 1998.

\bibitem{ck}
S.-J. Cheng and V.~G. Kac.
\newblock Conformal modules.
\newblock {\em Asian J. Math.}, 1(1):181--193, 1997.
\newblock Erratum: 2(1):153--156, 1998.

\bibitem{dk}
A.~D'Andrea and V.~G. Kac.
\newblock Structure theory of finite conformal algebras.
\newblock {\em Selecta Math. (N.S.)}, 4(3):377--418, 1998.

\bibitem{dlm_poiss}
C.~Dong, H.~Li, and G.~Mason.
\newblock Vertex {L}ie algebras, vertex {P}oisson algebras and vertex algebras.
\newblock In {\em Recent developments in infinite-dimensional Lie algebras and
  conformal field theory (Charlottesville, VA, 2000)}, volume 297 of {\em
  Contemp. Math.}, pages 69--96. Amer. Math. Soc., Providence, RI, 2002.

\bibitem{gd}
I.~M. Gel{\cprime}fand and I.~Ja. Dorfman.
\newblock Hamiltonian operators and algebraic structures associated with them.
\newblock {\em Funktsional. Anal. i Prilozhen.}, 13(4):13--30, 96, 1979.

\bibitem{jacobson_lie}
N.~Jacobson.
\newblock {\em Lie algebras}.
\newblock Dover Publications Inc., New York, 1979.
\newblock Republication of the 1962 original.

\bibitem{kac2}
V.~G. Kac.
\newblock {\em Vertex Algebras for Beginners}, volume~10 of {\em University
  Lecture Series}.
\newblock AMS, Providence, RI, second edition, 1998.

\bibitem{kac_fd}
V.~G. Kac.
\newblock Formal distribution algebras and conformal algebras.
\newblock In {\em XIIth International Congress of Mathematical Physics (ICMP
  '97) (Brisbane)}, pages 80--97. Internat. Press, Cambridge, MA, 1999.

\bibitem{primc}
M.~Primc.
\newblock Vertex algebras generated by {L}ie algebras.
\newblock {\em J. Pure Appl. Algebra}, 135(3):253--293, 1999.

\bibitem{freecv}
M.~Roitman.
\newblock On free conformal and vertex algebras.
\newblock {\em J. Algebra}, 217(2):496--527, 1999.

\bibitem{universal}
M.~Roitman.
\newblock Universal enveloping conformal algebras.
\newblock {\em Selecta Math. (N.S.)}, 6(3):319--345, 2000.

\bibitem{cfva}
M.~Roitman.
\newblock Combinatorics of free vertex algebras.
\newblock {\em J. Algebra}, 255(2):297--323, 2002.

\bibitem{zelmconf}
E.~Zelmanov.
\newblock On the structure of conformal algebras.
\newblock In {\em Combinatorial and computational algebra (Hong Kong, 1999)},
  volume 264 of {\em Contemp. Math.}, pages 139--153. Amer. Math. Soc.,
  Providence, RI, 2000.

\end{thebibliography}

\end{document}